\documentclass[11pt,leqno]{article}
\usepackage{amsmath}
\usepackage{amsfonts}
\usepackage{amssymb}
\usepackage{enumerate}
\usepackage{hyperref}
\title{Existence of weak solutions to stochastic evolution inclusions}
\author{Adam {\sc Jakubowski}\footnote{Nicholas Copernicus University, 
Faculty of Mathematics and Informatics, 
ul. Chopina 12/18, 
87--100 Toru\'n, Poland.  
E-Mail: adjakubo@mat.uni.torun.pl}\\
Mikhail {\sc Kamenski{\u \i}}\footnote{Departement of Mathematics, 
 State University of Voronezh, 
 Voronezh, Universitetskaja pl. 1, 
 394693, Russia. 
 E-Mail: mikhailkamenski@mail.ru}\\
Paul {\sc Raynaud de Fitte}\footnote{Laboratoire Rapha\"el Salem, UMR
   CNRS 6085, UFR Sciences , Universit\'e de Rouen, 76821 
   Mont Saint Aigan Cedex, France. E-Mail: prf@univ-rouen.fr}}
\date{August 6, 2004}

\renewcommand\paragraph[1]{\subsection*{#1}
                 \addcontentsline{toc}{subsection}{\numberline{}{#1}}
                           }
\newtheorem{theo}{Theorem}[section]
\newtheorem{prop}[theo]{Proposition}
\newtheorem{cor}[theo]{Corollary}
\newtheorem{lem}[theo]{Lemma}
\newtheorem{remark}[theo]{Remark}
\newenvironment{rem}{\begin{remark}\em}{\end{remark}}
\newtheorem{Bestimmung}[theo]{Definition}

\newtheorem{example}[theo]{Example}
\newtheorem{examples}[theo]{Examples}

\newtheorem{counterexample}[theo]{Counterexample}

\newtheorem{remdefi}[theo]{Remark and Definition}

\renewcommand\iff{if and only if\ }

\newcommand\preuv[1]{\par\medskip\noindent{\bf #1.}\ } 
\newcommand\preuve{\preuv{Proof}}
\newcommand\preuvof[1]{\preuv{Proof of {#1}}}
\newcommand\carre{\fbox{\rule{0em}{.3em}\rule{.3em}{0em}} \qquad}
\newcommand\fin{\hfill$\carre\qquad$\medskip\par} 
\newcommand\titre[1]{{\bf(#1)}}

\newcommand\eg{{\em {e.g.}}}
\newcommand\wlg{without loss of generality}

\newcommand\tq{;\,} 
\newcommand\foreach[1]{\CCO{\forall{#1}}\ }
\newcommand\thereis[1]{\CCO{\exists{#1}}\ }

\newcommand\ra{\rightarrow}
\newcommand\mt{\mapsto}
\newcommand\R{\mathbb R}
\newcommand\N{\mathbb N}
\newcommand\NN{\widehat{N}} 

\newcommand\restr[1]{ { \rule[-5pt]{0.4pt}{10pt}}_{#1} } 

\newenvironment{syst}{\left\{ \begin{array}{l}}{\end{array}\right.} 

\newcommand\abs[1]{\left\vert{#1}\right\vert}
\newcommand\accol[1]{\left\{{#1}\right\}}
\newcommand\CCO[1]{\left({#1}\right)}
\newcommand\croche[1]{\left[{#1}\right]}
\newcommand\norm[1]{\left\Vert{#1}\right\Vert}
\newcommand\un[1]{\,\rlap{{1}}\kern.22em \mbox{l}_{#1}} 

\newcommand\ctegeneral{K} 
\newcommand\grosscte{K} 

\newcommand\CteConv{C_{\text{\rm Conv}}}
\newcommand\CteCauchy[1]{M_{#1}} 

\newcommand\sgcontract{M} 
\newcommand\sgcoef{\beta}

\newcommand\growthconstant{C_{\text{\rm growth}}}

\newcommand\Lcte{k}

\newcommand\HyI{{\rm (HI)}} 
\newcommand\HySg{{\rm (HS)}}
\newcommand\HyF{{\rm (HFG)}}
\newcommand\HyFG{\HyF}

\newcommand\LL{L}
\newcommand\mapR{z} 

\newcommand\espgene{\ensuremath{\mathbb{E}}}
\newcommand\espmetricgeneric{\ensuremath{\mathbb{E}}} 
\newcommand\espmetricgene{\espmetricgeneric}
\newcommand\espmet{\espmetricgeneric}
\newcommand\esp{\ensuremath{\mathbb{H}}}

\newcommand\hilbw{\ensuremath{\mathbb{U}}}

\newcommand\Span{{\text{\rm Span}}} 
\newcommand\proj{P} 

\newcommand\HSz{{\mathbb L}}

\newcommand\compacts{\mathfrak{K}_c} 
\newcommand\compactspaconv{\mathfrak{K}}

\newcommand\Hausd[1][]{\mathop{\text{Hausd}}\nolimits_{#1}} 
\newcommand\dist[1][]{d_{#1}}

\newcommand\Times{[0,\Time]} 
\newcommand\Time{T}

\newcommand\bor[1]{{\mathcal B}_{#1}} 
\newcommand\lebesgue{\text{meas}} 

\newcommand\esprob{\Omega}                       
\newcommand\tribu{{\mathcal F}}                  
\newcommand\sstribu[1]{\tribu_{#1}}                    
\newcommand\pr{ \mathop{\text{\rm P}}\nolimits }          
\newcommand\oap{\ensuremath{(\esprob,\tribu, \pr)}}   

\newcommand\Ctribu{{\mathcal C}} 

\newcommand\laws[1]{{\mathcal M}^{1,+}\CCO{{#1}}}
\newcommand\Blop[1]{\text{\rm BL}_1({#1})} 
\newcommand\law[1]{ {\mathcal L}\CCO{#1} } 

\newcommand\comp{\mathcal{K}} 
\newcommand\comppp{\mathcal{R}} 
\newcommand\compC{\mathcal{K}} 
\newcommand\comppC{\mathcal{Q}} 
\newcommand\compppC{\mathcal{Z}} 
\newcommand\compppK{\mathfrak{K}} 

\newcommand\filtration{\ensuremath{(\tribu_t)}}

\newcommand\stbas{{\mathfrak F}}
\newcommand\prodspace{{\Times\times\esprob}}
\newcommand\predict{{\mathcal P}} 

\newcommand\expect{\mathop{\text{\bf E}}\nolimits}   

\newcommand\esprobb{\underline{\esprob}}
\newcommand\tribuu{\underline{\tribu}}

\newcommand\tas{{\mathfrak T}} 

\newcommand\skoroC{\text{\rm {C}}}
\newcommand\distC{d_{\infty}} 
\newcommand\Cb[1]{{\mathcal C}_b\CCO{#1}} 

\newcommand\Ell{\text{L}} 
\newcommand\randvar{\Ell^0} 
\newcommand\Ellp[1]{\Ell^{#1}}

\newcommand\bt{{\mathcal N}^p} 

\newcommand\BTC[2]{ \bt_c(\stbas,[0,#1];#2) } 
\newcommand\BTCC{ \bt_c(\underline{\stbas},[0,\Time];\esp) } 
\newcommand\BTCsimple[1]{\bt_c(#1)}

\newcommand\EProc{\Lambda} 
\newcommand\proc{X}
\newcommand\procc{\soll} 

\newcommand\ZZ{{\mathcal Z}} 

\newcommand\wien{W}

\newcommand\procHS{Z} 

\newcommand\smgr{S}
\newcommand\semigr[1]{\smgr(#1)} 
\newcommand\semigroupe{\CCO{\semigr{t}}_{t\geq 0}}

\newcommand\carathf{{\mathfrak f}}
\newcommand\carathg{{\mathfrak g}}
\newcommand\supercarath{{\mathfrak h}}

\newcommand\steiner{{\mathfrak S}}

\newcommand\youngs{{\mathcal Y}}
\def\dirac#1{{\delta}_{#1}} 
\def\ydirac#1{\underline{\delta}_{#1}} 

\newcommand\sol{X} 
\newcommand\ssol{\widetilde{\sol}} 
\newcommand\soll{Y}

\newcommand\Soll{\Xi}
\newcommand\initial{\xi} 

\newcommand\mncp{\mathop{\Psi}\nolimits} %
\newcommand\mncpp{\mathop{\Psi}\nolimits'}

\newenvironment{portemanteau}
  {\begin{enumerate}[{$1.$}]}{\end{enumerate}}
\newcommand{\qu}{\begin{portemanteau}}
\newcommand{\uq}{\end{portemanteau}}

\newenvironment{qromain}
  {\begin{enumerate}[{$(i)$}]}{\end{enumerate}}
\newcommand{\qrom}{\begin{qromain}}
\newcommand{\morq}{\end{qromain}}
\newcommand\qref[1]{(\ref{#1})}

\newenvironment{qualpha}
  {\begin{enumerate}[{$(a)$}]}{\end{enumerate}}
\newcommand{\qualph}{\begin{qualpha}}
\newcommand{\alphuq}{\end{qualpha}}

\begin{document}
\maketitle

\begin{abstract}
We prove the existence of a weak mild solution to the Cauchy problem 
for the semilinear stochastic differential inclusion in a Hilbert
space
$$dX_t\in AX_t\,dt+F(t,X_t)\,dt+G(t,X_t)\,dW_t$$
where 
$W$ is a Wiener process, 
$A$ is a linear operator which generates a $C_0$-semigroup, 
$F$ and $G$ are multifunctions with convex compact values satisfying
some growth condition and, with respect to the second variable, 
a condition 
weaker than the Lipschitz condition. 
The weak solution is constructed in the sense of Young
measures. 
\end{abstract}

\section{Introduction}
Ordinary and stochastic differential and inclusions in infinite dimensional
spaces have many important and interesting applications on which we
do not stop here, addressing the reader to the books and to the reviews (see for example 
\cite{ahmed91book,dapratozabczyk92book,bogachev95deterministic}). 

For ordinary differential equations,  after the famous
example of Dieudonn\'e 
\cite{dieudonne50exemples} and the fundamental results of A.~Godunov 
\cite{godunov72counter,godunov74peano}, 
it became evident that,  
in the case of infinite dimensional Banach spaces, 
it is necessary to suppose an auxilary
condition on the right hand side of the equation
\begin{equation}\label{eq:1}
x'=f(t,x)
\end{equation}
in order to have an existence theorem for the initial value problem with the initial
condition 
\begin{equation}\label{eq:2}
x(0)=x_0.
\end{equation}

It is well known that the estimation
\begin{equation}\label{eq:3}
\| f(t,x) - f(t,y)\|\le L(t, \|x-y\|),
\end{equation}       
where $L$ is a real function such that the integral inequality
\begin{equation*}
u(t)\le \int_0^t L(s, u(s))\,ds,
\end{equation*}
has a unique solution $u(t)=0,$ gives the existence condition mentioned below.
In \cite{kibenko-krasno-mamedov61one-sided} it was shown that it is 
possible to add to $f$ satisfying \eqref{eq:3} a continuous compact
operator. As generalisations of this fact many papers in the 
seventies 
were devoted to a condition of the form
\begin{equation}\label{eq:5}
\varphi (f(t, \Lambda))\le L(s,\varphi (\Lambda)), 
\end{equation}
where $\varphi $ is a measure of noncompactness (see e.g.~\cite{akprs92book}). The abstract
fixed point theorem for condensing operators (see \cite{akprs92book}) was successfully 
applied in this way.

In the end of the XXth century, in many papers (see
e.g.~\cite{ugowski85evolution,szufla86existence,banas87various}), 
it was shown 
that Condition \eqref{eq:5} implies the existence of solutions to
the semilinear equation
\begin{equation*}
x'=Ax+f(t,x)
\end{equation*}
and the semilinear inclusion (e.g.~\cite{obukhovski91controlled}, see
also the references in \cite{hu-papageorgiou00book})
\begin{equation}\label{eq:7}
x'\in Ax+f(t,x),
\end{equation}
with the initial condition \eqref{eq:2} and a linear operator $A$ generating a $C_0$-semigroup. 
In the inclusion \eqref{eq:7}, $f$ is a multivalued map with convex compact values.

In the same period from the 
seventies, 
it was remarked that the existence of a
weak solution 
and the condition of unicity of trajectories  implies
the existence of a strong solution 
to the stochastic differential equation
\begin{equation}\label{eq:8}
dX_t=a(t,X_t)\,dt+b(t,X_t)\,dW_t,
\end{equation}  
where $W$ is a standard Wiener process (here,``weak'' and ``strong'' are taken 
in the probabilistic sense; the existence of a weak solution 
with only continuous $a$ and $b$ was known since Skorokhod
\cite{skorokhod65book}). 
The direct
proof of the existence of a strong solution for \eqref{eq:8} 
using a convenient measure of noncompactness
when $a$ and $b$ satisfy a condition like \eqref{eq:3} was presented
in \cite{rodkina84heredity}, see also 
\cite{akprs92book}. 

The passage from the finite dimensional case to the infinite
dimentional case with $a$ and $b$ satisfying 
a Lipschitz condition is presented in
\cite{dapratozabczyk92book,bogachev95deterministic}  
see there the 
bibliography. 
The generalisation to the case of a semilinear stochastic differential equation
\begin{equation*}
dX_t=AX_t\,dt+f(t,X_t)dt+\sigma (t,X_t)\,dW_t,
\end{equation*}
deals with the estimations for the stochastic convolution operator
$$ \int_0^te^{(t-s)A}v(s)\,dW_s $$
(see
\cite{tubaro84burk,kotelenez84stopped,%
dapratozabczyk92book,dapratozabczyk92convolution,hausenblasseidler01stochconv}) 
and for Lipschitz $f$ and $\sigma$ is presented in
\cite{dapratozabczyk92book}. 

In the work of  
Da Prato and Frankowska \cite{dapratofrankowska94filippov},  
the existence result is proved for the semilinear stochastic inclusion
\begin{equation}\label{eq:10}
dX_t\in AX_t\,dt+f(t,X_t)\,dt+\sigma (t,X_t)dW_t,
\end{equation}
with multivalued $f$ and $\sigma$ which are Lipschitz with respect to
the Hausdorff metric.

In the present paper we aim to prove the existence of a ``weak''
solution (or a ``solution measure'') for \eqref{eq:10} in the case where
$f$ and $\sigma$ are compact valued multifunctions satisfying a  
condition like \eqref{eq:3} where the norm in the left hand side is
replaced by the Hausdorff metric. 
More precisely,  
we prove the existence of a ``weak'' mild  
solution $\sol$ to the Cauchy problem 
\begin{equation}\label{eq:generale}
\begin{syst}
d\sol_t\in A\sol_t+F(t,\sol_t)\,dt
+G(t,\sol_t)\,d\wien(t)\\
\sol(0)=\initial
\end{syst}
\end{equation}
where $\sol$ takes its values in a Hilbert space $\esp$, 
$\wien$ is a Brownian motion on a Hilbert space $\hilbw$, 
$A$ is a
linear operator on $\esp$ and $F$ and $G$ are 
multivalued mappings with compact
convex values, continuous in the second variable, 
and satisfy 
an assumption which is much more general than the usual Lispchitz one.

Instead of constructing the weak solution with the help of Skorokhod's
representation theorem, we define our weak solution as a Young
measure (ie as a {\em solution measure} in the sense of Jacod and
M\'emin \cite{jacod-memin81weakstrong}).

As in the case of the ordinary differential equation \eqref{eq:1} with the
right hand side satisfying 
\eqref{eq:5} (see \cite{kamenski72peano}) we apply the Tonelli scheme, 
and prove that its solutions are uniformly tight. 
We then 
pass to the limit and we obtain the existence of a weak mild solution
to the initial problem.

To get rid of the Lipschitz assumption on $F$ and $G$ has a cost: 
not only do we obtain ``weak'' solutions, but  
our techniques based on compactness lead us to consider mappings $F$
and $G$ with convex compact values, 
whereas the more geometric methods of  
\cite{dapratofrankowska94filippov} deal  
with unbounded closed valued mappings.     
Furthermore, in our work, the multifunctions $F$ and $G$ 
are deterministic, whereas they are random in
\cite{dapratofrankowska94filippov}. 

\section{Formulation of the problem, statement of the result}
\paragraph{Notations}

Throughout, 
$0<\Time<+\infty$ is a fixed time 
and 
$\stbas=(\esprob,\tribu, \filtration_{t\in\Times},\pr)$ is a 
stochastic basis satisfying the usual conditions.  
The $\sigma$--algebra of predictable subsets of $\prodspace$ is denoted by
$\predict$.

If $\espmetricgeneric$ is a separable metric space, we denote 
$\bor{\espmetricgeneric}$ its Borel $\sigma$--algebra, and by
$\laws{\espmetricgeneric}$ the space of probability laws on
$(\espmetricgeneric,\bor{\espmetricgeneric})$, endowed with the usual narrow (or weak)
topology. 
The law of a random element $X$ of $\espmetricgeneric$ is denoted by
$\law{X}$. 
The space of random elements of $\espmetricgeneric$ defined on
$(\esprob,\tribu)$ is denoted by $\randvar(\esprob;\espmetricgeneric)$
or simply $\randvar(\espmetricgeneric)$, and
it is endowed with the topology of convergence in
$\pr$--probability. We identify random elements which are equal
$\pr$--a.e. 
Recall that a subset $\Lambda$ of $\randvar(\espmetricgeneric)$ is
said to be tight 
if, 
for any $\delta>0$, there exists a compact subset
$\mathcal K$ of $\espmetricgeneric$ such that, for every
$X\in\Lambda$, $\pr(X\in {\mathcal K})\geq 1-\delta$.

\paragraph{Spaces of processes}
In all this paper, $p> 2$ is a fixed number. 
For any $t\in\Times$, 
we denote by $\BTC{t}{\esp}$ the space of continuous
$\filtration$--adapted $\esp$--valued
processes $\sol$ such that 
\begin{equation*}
\norm{\sol}^p_{ \BTC{t}{\esp} }
:=\expect\CCO{\sup_{0\leq s\leq t}\norm{\sol(s)}_\esp^p}<+\infty.
\end{equation*}

For any $t\in\Times$ 
we denote by 
$\skoroC([0,t];\esp)$ the space of continuous mappings from $[0,t]$
to $\esp$. The space $\skoroC([0,t];\esp)$ is endowed with the
topology of uniform convergence defined by the distance 
$\distC(u,v)=\sup_{s\in[0,t]}\norm{u(s)-v(s)}$. 
The space $\BTC{t}{\esp}$ is thus a closed separable subspace of the space
$\randvar\CCO{\skoroC([0,t];\esp)}$ of $\skoroC([0,t];\esp)$--valued
random variables. 

\paragraph{Measures of noncompactness}
If $(\espmetricgene,d)$ is a metric space and
$\EProc\subset\espmetricgene$, we say that a subset $\EProc'$ of
$\espmetricgene$ is an {\em$\epsilon$--net} of $\EProc$ if 
\begin{equation*}
  \inf_{\proc\in\EProc}d(\proc,\EProc')\leq \epsilon 
\end{equation*}
(note that $\EProc'$ is not necessarily a subset of $\EProc$). 
Let $\EProc$ be a subset of $\BTC{\Time}{\esp}$. 
For any $s\in[0,\Time]$, we denote 
by $\EProc\restr{s}$ the subset of $\BTC{s}{\esp}$ of
restrictions to $[0,s]$ of elements of $\EProc$. 
We denote 
\begin{equation*}
   \mncp(\EProc)(s):=\inf\accol{\epsilon>0\tq
   \text{$\EProc\restr{s}$ has a tight $\epsilon$--net in $\BTC{s}{\esp}$}}. 
\end{equation*}
We shall see in Lemma \ref{lem:G}
that $\EProc$ is tight \iff $\mncp(\EProc)(\Time)=0$. 
Note also that the mapping $s\mt\mncp(\EProc)(s)$ is nondecreasing. 
We denote 
\begin{equation*}
  \mncp(\EProc):=\CCO{\mncp(\proc)(s)}_{0\leq s\leq \Time}.
\end{equation*}
The family $\CCO{\mncp(\proc)(s)}_{0\leq s\leq t}$ is called the
{\em measure of noncompactness} of $\EProc$ (see
\cite{akprs92book,koz01book} about measures of noncompactness).

\paragraph{Spaces of closed compact subsets}
For any metric space $\espgene$, the set of nonempty compact 
(resp.~nonempty compact convex) 
subsets of $\espgene$ is denoted by $\compactspaconv(\espgene)$ (resp.~$\compacts(\espgene)$). 
We endow $\compactspaconv(\espgene)$ and its subspace
$\compacts(\espgene)$ with the Hausdorff distance 
$$\Hausd[\espgene](B,C):=\max\CCO{\inf_{b\in B}d(b,C),\inf_{c\in C}d(c,B)}.$$ 
Recall that, if $\espgene$ is Polish, then $\compacts(\espgene)$ is
Polish too (see \eg~\cite[Corollary II-9]{castaingvaladier77book}). 


\paragraph{Stable convergence and Young measures}
Let $\espmetricgeneric$ be a complete metric space.  
We denote by $\youngs(\esprob,\tribu,\pr;\espmetricgeneric)$ (or simply
$\youngs(\espmetricgeneric)$) the set of measurable mappings 
$$\mu :\,\left\{ 
\begin{matrix}\esprob&\mt&\laws{\espmetricgeneric}\\
\omega&\mt&\mu_\omega
\end{matrix}
\right.$$ 
Each element $\mu$ of $\youngs(\esprob,\tribu,\pr;\espmetricgeneric)$ can be
identified with the measure $\widetilde{\mu}$ on 
$(\esprob\times\espmetricgeneric,\tribu\otimes
\bor{\espmetricgeneric})$ 
defined by $\widetilde{\mu}(A\times
B)=\int_A\mu_\omega(B)\,d\pr(\omega)$ (and the mapping
$\mu\mt\widetilde{\mu}$ is onto, 
see \eg~\cite{valadier73desi}). In the sequel, we shall use freely this
identification. For instance,  if $f:\,\esprob\times\espmetricgeneric\ra\R$ 
is a bounded measurable mapping, the notation $\mu(f)$ denotes 
$\int_\esprob \mu_\omega(f(\omega,.))\,d\pr(\omega)$. 
The elements of $\youngs(\esprob,\tribu,\pr;\espmetricgeneric)$ are called 
{\sl Young measures} on $\esprob\times\espmetricgeneric$. 

Let $\Cb{\espmetricgeneric}$ denote the set of continuous bounded real
valued functions defined on $\espmetricgeneric$. 
The set $\youngs(\esprob,\tribu,\pr;\espmetricgeneric)$ is endowed
with a metrizable topology, such that 
a sequence $(\mu^n)$ of Young measures converges to a Young measure $\mu$ if, 
for each $A\in\tribu$ and  each $f\in\Cb{\espmetricgeneric}$, 
the sequence $(\mu^n(\un{A}\otimes f))$ converges to
$\mu(\un{A}\otimes f)$ (where $\un{A}$ is the indicator function of $A$). 
We then say that $(\mu^n)$ converges {\em stably} to $\mu$. 

Each element $X$ of $\randvar(\esprob;\espmetricgeneric)$ can (and
will sometimes) be
identified with the Young
measure 
$$\ydirac{X} :\,\omega\mt\dirac{X(\omega)},$$ 
where, for any $x\in\espmetricgeneric$, $\dirac{x}$
denotes the probability concentrated on $x$. 
Note that the restriction to  $\randvar(\esprob;\espmetricgeneric)$ of
the topology of stable convergence is the topology of convergence in
probability. 
If $(X_n)$ is a tight sequence in
$\randvar(\esprob;\espmetricgeneric)$,  
then, by Prohorov's compactness criterion 
for Young measures \cite{balder89Prohorov,cc-prf-valadier04book},
each subsequence of $(X_n)$ has a further subsequence, say $(X'_n)$,  
which converges stably to some $\mu\in\youngs(\espmetricgeneric)$, 
that is, for every $A\in\tribu$ and every $f\in\Cb{\espmetricgeneric}$, 
$$\lim_n\int_A f(X'_n)\,d\pr=\mu(\un{A}\otimes f).$$
This entails in particular that $(X'_n)$ converges in ditribution to
the measure $\mu(\esprob\times.)\in\laws{\espmetricgeneric}$.

See \cite{valadier94course} or \cite{balder00lectures} for 
 an introduction to Young measures and their applications.

\paragraph{Hypothesis}
In the sequel, 
we are given two separable Hilbert spaces $\esp$ and $\hilbw$ and 
an $\filtration_{t\in\Times}$--Brownian
motion $\wien$ (possibly cylindrical) on $\hilbw$. 
We denote by $\HSz$ the space of Hilbert--Schmidt 
operators from $\hilbw$ to $\esp$.

We shall consider the following hypothesis: 

\begin{list}{}{\leftmargin 4em\labelwidth 2em}
\item[\HySg] $A$ is the generator of a $\mathcal{C}_0$ semigroup
  $\semigroupe$. 
  In particular 
  (see \eg~\cite[Theorem 1.3.1]{ahmed91book}), there exist
  $\sgcontract>0$ and $\sgcoef\in]-\infty,+\infty[$ such that, for
  every $t\geq 0$, 
  \begin{equation*}
  \norm{\semigr{t}}\leq\sgcontract e^{\sgcoef t}.
  \end{equation*}  
  For $t\in[0,\Time]$, we denote $\CteCauchy{t}=
                           \sup_{0\leq s\leq t}\sgcontract e^{\sgcoef s}$.

\item[\HyF] 
  $F :\,\esprob\times\Times\times\esp\ra \compacts(\esp)$
  and $G :\,\esprob\times\Times\times\esp\ra \compacts(\HSz)$ are
  measurable mappings which satisfy
  the following conditions: 
\qrom
\item \label{hyp:growth}
  There exists a constant $\growthconstant>0$ such that, for all
  $(t,x)\in\Times\times\esp$, 
  \begin{align*}
    \Hausd[\esp](0,F(t,x))&\leq \growthconstant(1+\norm{x})\\
    \Hausd[\HSz](0,G(t,x))&\leq \growthconstant(1+\norm{x}).
  \end{align*}

\item\label{hyp:FL}  
  For all $(t,x,y)\in\Times\times\esp\times\esp$,  
  \begin{align*}
   \CCO{\Hausd[\esp](F(t,x),F(t,y))}^p
               &\leq \LL(t,\norm{x-y}^p)\\
  \CCO{\Hausd[\HSz](G(t,x),G(t,y))}^p
              &\leq \LL(t,\norm{x-y}^p),
  \end{align*}
where $\LL :\,[0,\Time]\times [0,+\infty]\ra
[0,+\infty]$ is a given continuous mapping such that 
\qualph 
\item  for every $t\in[0,\Time]$, the mapping $\LL(t,.)$ is
  nondecreasing and convex,
\item\label{hyp:L} 
  for every measurable mapping $\mapR :\,
  [0,\Time]\ra[0,+\infty] $ and for every constant $\ctegeneral>0$, 
  the following implication holds true:
  \begin{equation}\label{eq:hyp:L}
    \croche{\foreach{t\in[0,\Time]}
           \mapR(t)\leq \ctegeneral \int_0^t \LL(s,\mapR(s))\,ds 
           }\Rightarrow \mapR =0.
  \end{equation}
In particular, we have $\LL(0)=0$, thus Hypothesis \HyFG-\qref{hyp:FL}
entails that, for each $t\in\Times$, the mappings $F(t,.)$ and  
$G(t,.)$ are continuous for the Hausdorff distances $\Hausd[\esp]$ and
$\Hausd[\HSz]$ respectively. 
Such a function $\LL$ is considered in
e.g.~\cite{yamada81successive,rodkina84heredity,manthey90convergence,akprs92book,%
taniguchi92successive,barbu98local-global,barbu-bocsan02approx}. 
Concrete examples can be found in \cite[Section 6 of Chapter 3]{hartman73book}. 

\alphuq
\morq

\item[\HyI]
  $\initial\in\Ellp{p}(\esprob,\sstribu{0},\pr\restr{\sstribu{0}};\esp)$.  

\end{list}
Recall that, under Hypothesis \HySg, 
there exists a constant $\CteConv$ such that, for any predictable process
  $\procHS\in\Ellp{p}(\esprob\times\Times;\HSz)$, we have  
\begin{equation}\label{stochconvcondense}
\expect\croche{\sup_{s\leq t}\norm{ 
\int_0^s \semigr{s-r}\procHS(r)\,d\wien(r)  }^p}
\leq \CteConv\,t^{(p/2)-1}\,\expect\int_0^t\norm{\procHS(s)}_{\HSz}^p\,ds 
\end{equation}
(see 
\cite{dapratozabczyk92convolution,dapratozabczyk92book}, and also 
\cite{hausenblasseidler01stochconv} for a strikingly short proof in the case
when $\smgr$ is contractive).  

\paragraph{Weak and strong mild solutions}
We say that $\sol\in\BTC{\Time}{\esp}$  is a  {\em (strong) mild solution}
to Equation \eqref{eq:generale} if 
there exist two predictable processes $f$ and $g$ defined on $\stbas$ satisfying 
\begin{equation}\label{eq:ids}
\left\{
\begin{aligned}
&\sol(t)= \semigr{t}\initial+\int_0^t
\semigr{t-s}f(s)\,ds+\int_0^t 
\semigr{t-s} g(s)\,d\wien(s)\\
&f(s)\in F(s,\sol(s))\ \pr\text{-a.e.}\\
&g(s)\in G(s,\sol(s))\ \pr\text{-a.e.}. 
\end{aligned}\right.
\end{equation}
So, ``mild solution'' refers to the variation of constant formula,
whereas 
``strong'' refers to the fact that the solution is defined on the
given stochastic basis. 

We say that a process $\sol$ is a {\em weak mild solution} 
or a {\em mild solution-measure} 
to 
\eqref{eq:generale} if there exists a stochastic basis 
$\underline{\stbas}$\allowbreak
$=(\esprobb,$\allowbreak$\tribuu,$\allowbreak$(\tribuu_t)_t,$\allowbreak$\mu)$
satisfying the following conditions:  
\begin{enumerate}
\item 
$\esprobb$ has the form $\esprobb=\esprob\times\esprob'$, 
$\tribuu=\tribu\otimes\tribu'$ for some $\sigma$--algebra $\tribu'$ on
$\esprob'$, $\tribuu_t=\tribu_t\otimes\tribu'_t$ for some right
continuous filtration  $(\tribu'_t)$ on $(\esprob',\tribu')$, and the
probability $\mu$ satisfies $\mu(A\times\esprob')=\pr(A)$ for every
$A\in\tribu$. 

\item The process $\wien$ is a Brownian motion on $\underline{\stbas}$ 
  (we identify here every random variable $X$ on $\Omega$ with the
  random variable $(\omega,\omega')\mt X(\omega)$ defined on
  $\esprobb$).

\item
  $X\in\BTCC$ 
  and there exist two predictable processes $f$ and $g$ defined on $\underline{\stbas}$ 
  satisfying \eqref{eq:ids}. 
  
\end{enumerate}
  The terminology {\em solution-measure} is that of 
  \cite{jacod-memin81weakstrong}. 
  If $\tribu'$ is the Borel $\sigma$--algebra of some
  topology on $\esprob'$, a solution-measure can also be seen as a
  Young measure. 
  This is the point of view adopted by Pellaumail
  \cite{pellaumail81solfaibles,pellaumail81weak}, who calls Young
  measures {\em rules}. 

\paragraph{Main result and corollaries}
We can now state the main result of this paper. 
The proofs will be given in Section \ref{sect:proofs}. 
\begin{theo}\label{theo:main}
\titre{Main result} 
Under Hypothesis \HySg, \HyFG\ and \HyI, Equation \eqref{eq:generale} has a
weak mild solution. 
\end{theo}
An easy adaptation of our reasoning also yields, as a by-product, a
well-known strong existence
result: 
\begin{prop}\label{prop:single-valued}
\titre{Strong existence in the single valued case \cite{barbu98local-global,barbu-bocsan02approx}}
Under Hypothesis \HySg, \HyFG\ and \HyI,
 if furthermore $F$ and $G$ are single-valued, then \eqref{eq:generale} has a
 strong mild solution.
\end{prop}
Using the Steiner point for the choice of selections of $F$ and $G$,
we deduce the following 
\begin{prop}\label{cor:finitedim}
\titre{Strong existence in the finite dimensional case}
Under Hypothesis \HySg, \HyFG\ and \HyI,
 if furthermore
$\esp$ and $\hilbw$ are finite dimensional, then \eqref{eq:generale} has a
 strong mild solution. 
\end{prop}


\section{Preliminary results} 

\paragraph{Tightness results and boundedness results}
We start with a very simple and useful lemma. 
\begin{lem}\label{lem:G}
\titre{A tightness criterion}
Let $(\espmet,\dist)$ be a separable complete metric space. 
Let $\EProc$ be a set of random
elements of $\espmet$ defined on $\oap$. 
Let $r\geq 1$. 
Assume that, for every $\epsilon>0$, there exists a tight subset
$\EProc_\epsilon$ such that
\begin{equation*}
\EProc\subset\EProc_\epsilon^{[\epsilon]}:=\{X\in \randvar\CCO{\espmet}\tq \exists
Y\in\EProc_\epsilon,\ \expect\dist^r(X,Y)<\epsilon \}.
\end{equation*} 
Then $\EProc$ is tight. 
\end{lem}
\preuve 
From Jensen inequality, we only need to prove Lemma
\ref{lem:G} for $p=1$. 
Indeed, we have $\expect\dist(X,Y)\leq (\expect\dist^r(X,Y))^{1/r}$
thus 
$$\EProc\subset\{X\in \randvar{(\espmet)}\tq \exists
Y\in\EProc_\epsilon,\ \expect\dist(X,Y)<\epsilon^{1/r}\}.$$

Let $\Blop{\dist}$ be the set of all mappings $f :\,\espmet\ra[0,1]$ wich are
1--Lischitz with respect to $\dist$. 
Let $\beta$ denote the Dudley distance on $\laws{\espmet}$, 
that is, for all $\mu,\nu\in\laws{\espmet}$,
\begin{equation*}
\beta(\mu,\nu)=\sup_{f\in\Blop{\dist}}\mu(f)-\nu(f). 
\end{equation*}
It is well known that the narrow topology on  $\laws{\espmet}$ is induced
by $\beta$ and that $\beta$ is complete, see \eg~\cite{dudley02book}. 
For every $X\in\EProc$, there exists $Y\in\EProc_\epsilon$ such that
$\expect\dist(X,Y)<\epsilon$, which implies 
$\beta(\law{X},\law{Y})<\epsilon$. 
We thus have  
\begin{equation*}
\accol{\law{X}\tq X\in\EProc}
\subset \cap_{\epsilon>0}\accol{\law{X}\tq X\in\EProc_\epsilon}^{[\epsilon]},
\end{equation*}
where, for any $\Xi\subset\laws{\espmet}$ and any $\epsilon>0$,
$\Xi^{[\epsilon]}=\{\mu\in\laws{\espmet}\tq \exists \nu\in\Xi,\
\beta(\mu,\nu)<\epsilon\}$.  
This proves that $\{\laws{X}\tq X\in\EProc\}$ is totally bounded for
$\beta$, thus relatively compact in the narrow topology. 
\fin

 \begin{lem}\label{lem:Zero}
  Let $\EProc$ be a set of continuous adapted processes on
  $\esp$. Assume that each element of $\EProc$ is in
  $\Ellp{p}(\esprob\times\Times;\esp)$ and 
that 
 that $\EProc$, considered
 as a set of $\skoroC(\Times;\esp)$--valued random variables, is tight.  
 Assume furthermore that $G$ satisfies Hypothesis
 \HyFG-\qref{hyp:growth}. 
 Let $\Soll$ be the set of processes $\soll$ of the form 
 $$\soll(t)=\int_0^t\semigr{t-s}g(s)\,d\wien(s)$$
 where $g$ is $\predict$--measurable and 
 $g(s)\in G(s,\proc(s))$ a.e.~for some $\proc\in\EProc$. 
 The set $\Soll$ is a tight set of $\skoroC(\Times;\esp)$--valued random
 variables. 
\end{lem}
\preuve
We will prove Lemma \ref{lem:Zero} through a series of reductions.

{\em First step:} {\em We can assume \wlg\ that there exists a compact
subset $\comppC$ of $\skoroC(\Times;\esp)$ such that, for each
$\proc\in\EProc$ and for each $\omega\in\esprob$,
$\proc(\omega,.)\in\comppC$.}  

Indeed, assume that Lemma \ref{lem:Zero} is true under this additional hypothesis. 
Let $\epsilon>0$. There exists a compact subset
$\comppC_\epsilon$ of $\skoroC(\Times;\esp)$ such that, for each
$\proc\in\EProc$, $\pr(\proc\in\comppC_\epsilon)\geq 1-\epsilon$. 
For every $\proc\in\EProc$, there exists a measurable subset
$\esprob_\proc$ of $\esprob$ such that $\pr(\esprob_\proc)\geq
1-\epsilon$ and, for every $\omega\in\esprob_\proc$,
$\proc(\omega)\in\comppC_\epsilon$. 
For each $\proc\in\EProc$ and each $t\in\Times$,  let us denote 
$$\proc^\epsilon(t)=\begin{cases}
\proc(t\wedge \tau)&\text{ if }\tau<+\infty\\
\proc(0)&\text{ if }\tau=+\infty,
\end{cases}
$$
where $\tau(\omega)$ is the infimum of all $s\in[0,\Time]$ such that 
there exists $u\in\comppC$ which coincides with $\proc(\omega,.)$ on
$[0,s]$  
 (we take $\inf\emptyset=+\infty$). 
The process $\proc^\epsilon$ is continuous and adapted, thus predictable. 
Furthermore, 
for every $\omega\in\esprob_\proc$, we
have $\proc(\omega,.)=\proc^\epsilon(\omega,.)$. 

Let $\EProc^\epsilon=\{\proc^\epsilon\tq\proc\in\EProc\}$. 
The set $\EProc^\epsilon$ is tight, 
thus, from our hypothesis, the set 
$\Soll^\epsilon$, obtained by replacing $\EProc$ by $\EProc^\epsilon$
in the definition of $\Soll$, is tight. 
There exists a compact subset $\compppC_\epsilon$ of $\skoroC(\Times;\esp)$ such that, for each
$\procc^\epsilon\in\Soll^\epsilon$, 
$\pr(\procc^\epsilon\in\compppC_\epsilon)\geq 1-\epsilon$. 
Let $\procc\in\Soll$. The process $\procc$ has the form 
$$\procc(t)=\int_0^t\semigr{t-s}g(s)\,d\wien(s)$$
 where $g$ is $\predict$--measurable and 
$g(s)\in G(s,\proc(s))$ a.e.~for some
 $\proc\in\EProc$. 
The set $\{\proc\not=\proc^\epsilon\}$ is predictable, thus there
 exists a predictable process $g^\epsilon$ such that 
 $g^\epsilon(s)\in G(s,\proc^\epsilon(s))$
 a.e.~ and $g^\epsilon(t)=g(t)$
 for $t\leq \tau$ (we can construct $g^\epsilon$ as a selection of the
 predictable multifunction $H$ defined by
 $H(t)=\{g(t)\}$ if $t\leq\tau$ and
 $H(t)=G(t,\proc^\epsilon(t))$ otherwise). 
Let $\procc^\epsilon\in\Soll^\epsilon$ be defined by 
$$\procc^\epsilon(t)=\int_0^t\semigr{t-s}g^\epsilon(s)\,d\wien(s).$$ 
We have 
\begin{align*}
\procc(t)-\procc^\epsilon(t)
&=\int_{t\wedge\tau}^t \semigr{t-s}(g(s)-g^\epsilon(s))\,d\wien(s)
\end{align*}
We thus have
$\procc=\procc^\epsilon$ on $\esprob_\proc$. 
This shows that, for every $\procc\in\Soll$, we have 
\begin{align*}
\pr(\procc\not\in\compppC_\epsilon)
&=\pr(\procc=\procc^\epsilon\text{ and
}\procc^\epsilon\not\in\compppC_\epsilon)+
\pr(\procc\not=\procc^\epsilon\text{ and
}\procc\not\in\compppC_\epsilon)\\
&\leq
\pr(\procc^\epsilon\not\in\compppC_\epsilon)+\pr(\procc^\epsilon\not=\procc)
\leq 2\epsilon. 
\end{align*}
Thus $\Soll$ is tight. 

\medskip
{\em Second step: We can furthermore assume \wlg\ that there exists a
compact subset $\comp$ of $\HSz$ such that
$G(t,\proc(\omega,t))\subset\compC$ a.e.~for all $\proc\in\EProc$.} 

Assume that Lemma \ref{lem:Zero} holds under this hypothesis. 

For any adapted process $g\in\Ellp{p}(\esprob\times\Times;\HSz)$, 
we denote by $\soll^{(g)}$ the process
defined by 
$$\soll^{(g)}(t)=\int_0^t\semigr{t-s}g(\omega,s)\,d\wien(s). $$ 
We denote by $\Gamma$ the set of 
$\predict$--measurable $\HSz$--valued processes $g$ such that  
$g(\omega,t)\in G(s,\proc(t))$ for every $t\in\Times$ a.e.~for some
 $\proc\in\EProc$.

Let $\comppp_0$ be a compact subset of $\esp$ such that
$\comppC(t)\subset \comppp_0$ for all $t\in\Times$ 
(where $\comppC$ is as in Step 1 and 
$\comppC(t)=\{u(t)\tq u\in\comppC\}$). 
The multifunction 
$t\mt G(t,\comppp_0)$ is measurable and has compact values in $\HSz$. 
Let $\epsilon>0$. 
There exists a compact subset ${\compppK}_\epsilon$ of
$\compacts(\HSz)$ and a measurable subset $I_\epsilon$ of $\Times$
such that, denoting by $\lebesgue$ the Lebesgue measure on $\Times$, 
\begin{gather}
\lebesgue\CCO{\Times\setminus I_\epsilon}\leq \epsilon,
\quad\text{ and }\quad
\foreach{t\in I_\epsilon}
G(t,\comppp_0)\in{\compppK}_\epsilon. \label{eq:Iepsilon}
\end{gather} 
Now, from a well known characterization of the compact subsets of 
$\compactspaconv(\HSz)$
(\cite[Theorem 2.5.2]{michael51topol}, see also  
\cite[Theorem {3.1}]{christensen74book} for the converse 
implication), the set 
$$\comp_\epsilon=\bigcup_{K\in{\compppK}_\epsilon}K$$ 
is a compact subset of $\HSz$. 
Let us define a multifunction 
$$G_\epsilon :\,\,\left\{ 
\begin{array}{lcl}
\Times\times\esp&\ra&\compacts(\HSz)\\
(t,x)&\mt&\un{I_\epsilon}(t)G(t,x).
\end{array}
\right.$$
For every 
$g\in\Gamma$ and every $t\in\Times$, let us set 
$$
g_\epsilon(t)=\un{I_\epsilon}(t)g(t). 
$$
The process $g_\epsilon$ is $\predict$--measurable and satisfies,  for some
 $\proc\in\EProc$, 
 \begin{equation*}
 g_\epsilon(t)\in G_\epsilon(t,\proc(t))\subset\comp_\epsilon \ 
\text{ a.e.~for every $t\in\Times$}
 \end{equation*}
 (we assume \wlg\ that $0\in\comppp_0$). 
Thus, from our hypothesis, the set 
\begin{gather*}
\Soll_\epsilon=\{\soll^{g_\epsilon}\tq g\in \Gamma \} 
\end{gather*} 
is tight. Let $R_0=\sup_{x\in\comppp_0}\norm{x}_\esp$. 
By \HyFG-\qref{hyp:growth}, we have, for every $g\in\Gamma$ and every $t\in\Times$, 
\begin{equation*}
\norm{g(t)-g_\epsilon(t)}\,\begin{cases}
=0&\text{if }t\in I_\epsilon\\
=\norm{g(t)}\leq 
\growthconstant\,(1+R_0)&\text{if }t\not\in I_\epsilon.
\end{cases}
\end{equation*}
Using \eqref{stochconvcondense}, we thus have
\begin{align*}
\Hausd[\BTC{\Time}{\esp}]\CCO{\Soll,\Soll_\epsilon}
&\leq \sup_{g\in\Gamma}\,\CteConv\,\Time^{(p/2)-1}\, 
       \expect\int_0^t\norm{g(s)-g_\epsilon(s)}_{\HSz}^p\,ds \\
&\leq \CteConv\,\Time^{(p/2)-1}
  \int_0^t \bigl(\growthconstant\,(1+R_0)\un{\Times\setminus
    I_\epsilon(s)}\bigr)^p\,ds\\
&\leq \CteConv\,\Time^{(p/2)-1}\,\growthconstant^p\,(1+R_0)^p\,\epsilon^p. 
\end{align*}
From Lemma \ref{lem:G}, we deduce that $\Soll$ is tight.

\medskip
{\em Third step: We can also assume \wlg\ that $A$ is a bounded
operator on $\esp$.} 

Let $A_n$ ($n>\sgcoef$) be the Yosida approximations of $A$. 
We are going to prove that 
\begin{equation}\label{eq:step3}
\sup_{g\in \Gamma}\expect\sup_{0\leq t\leq \Time}
\norm{
\int_0^t \CCO{\semigr{t-s}-e^{(t-s)A_n}}g(s)\,d\wien(s)  
}^p\ra 0\ \text{ as }n\ra+\infty.
\end{equation}

Let $D(A)$ be the domain of $A$. We have $\overline{D(A)}=\esp$ thus, 
for each $\epsilon>0$, we can find a
finite $\epsilon$--net $\comp_\epsilon$ of $\comp$ which lies in
$D(A)$. Then, for each $g\in\Gamma$, we can define a predictable
$\comp_\epsilon$--valued process $g_\epsilon$ such that
$\norm{g-g_\epsilon}_\infty\leq \epsilon$. 
For each $n$, the semigroup $e^{tA_n}$ satisfies an inequality similar
to \eqref{stochconvcondense} with same constant $\CteConv$, because
$\CteConv$ depends only on the parameters $\sgcontract$ and
$\sgcoef$ in Hypothesis \HySg. We thus have 
\begin{align*}
\sup_{g\in \Gamma}\expect\sup_{0\leq t\leq \Time}
\norm{
\int_0^t \semigr{t-s}\CCO{g(s)-g_\epsilon(s)}\,d\wien(s)  
}^p
&\leq \CteConv\Time^{p/2}\,\epsilon^p\\
\sup_{g\in \Gamma}\expect\sup_{0\leq t\leq \Time}
\norm{
\int_0^t e^{(t-s)A_n}\CCO{g(s)-g_\epsilon(s)}\,d\wien(s)  
}^p&\leq \CteConv\Time^{p/2}\,\epsilon^p.
\end{align*}
Therefore, we only need to prove \eqref{eq:step3} in the case when
$\comp\subset D(A)$.

For every $x\in D(A)$ and every integer $n>\sgcoef$, we have 
\begin{align*}
\CCO{\semigr{t-s}-e^{(t-s)A_n}}x
&=\croche{\semigr{t-s-r}\,e^{r
    A_n}}_0^{t-s}\,x
=\int_0^{t-s}
\semigr{t-s-r}\,e^{\tau A_n}\CCO{A_n-A}x\,dr.
\end{align*}
 Thus, assuming that $\comp\subset D(A)$, and 
 denoting by $C$ a constant which may difer from line to line, 
 we have, for every
 $g\in\Gamma$, using the stochastic Fubini theorem (see
 \cite{dapratozabczyk92book}) and the convolution inequality
 \eqref{stochconvcondense}, 
\begin{multline*}
\expect\sup_{0\leq t\leq \Time}
\norm{
\int_0^t \CCO{\semigr{t-s}-e^{(t-s)A_n}}g(s)\,d\wien(s)  
}^p\\
\begin{aligned}
&=\expect\sup_{0\leq t\leq \Time}
\norm{
\int_0^t \int_0^{t-s} \semigr{t-s-r}\,e^{r A_n}\CCO{A_n-A}g(s)\,dr\,d\wien(s)  
}^p\\
&=\expect\sup_{0\leq t\leq \Time}
\norm{
\int_0^t e^{r A_n}\,\int_0^{t-r} \semigr{t-s-r}\,e^{r A_n}\CCO{A_n-A}g(s)\,d\wien(s) \,dr 
}^p\\
&\leq C\expect\sup_{0\leq t\leq \Time}
\int_0^t \norm{
\int_0^{t-r} \semigr{t-s-r}\,\CCO{A_n-A}g(s)\,d\wien(s) 
}^p dr \\
&\leq C\int_0^\Time\expect\sup_{0\leq t\leq \Time}
\norm{
\int_0^{t-r} \semigr{t-s-r}\,\CCO{A_n-A}g(s)\,d\wien(s) 
}^p dr \\
&\leq C\int_0^\Time\int_0^\Time\expect
\norm{
\CCO{A_n-A}g(s)}^p
ds\,dr.
\end{aligned} 
\end{multline*}
But, for every $x\in D(A)$, we have $(A_n-A)x\ra 0$. So,
using the compactness of $\comp$ and 
 Lebesgue's dominated convergence theorem, we obtain \eqref{eq:step3}.

From Lemma \ref{lem:G}, we conclude that we only need to check 
Lemma \ref{lem:Zero} for the semigroups $e^{tA_n}$ ($n>\sgcoef$), which amounts to
check Lemma \ref{lem:Zero} in the
case when $D(A)=\esp$. 
In this case, $S(t)$ is the exponential $e^{tA}$ in the
usual sense. 

\medskip
{\em Fourth step: We can assume \wlg\ that $\esp$ is finite
dimensional.} 

Let $(e_n)$ be an orthonormal basis of $\esp$. For each $n$, let 
$\esp_n=\Span\CCO{e_1,\dots,e_n}$ and let $\proj_n$ be the orthogonal 
projection from $\esp$ onto $\esp_n$. 
Let $\Gamma$ be any contour around the spectrum of $A$, say $\Gamma$
is a circle $C(0,\rho)$. 
Denoting by $R$ the resolvent operator, we have, for any $g\in\Gamma$,
\begin{multline*}
\int_0^t \CCO{e^{(t-s)A}-e^{(t-s)P_nA}}P_ng(s)\,d\wien(s)  
\\
\begin{aligned}
&=\int_0^t\frac{1}{2\pi i}\int_\Gamma
e^{\lambda(t-s)}\,\CCO{R(\lambda,A)-R(\lambda,P_nA)}\,d\lambda\,
P_ng(s)\,d\wien(s)\\
&=\frac{1}{2\pi i}\int_\Gamma\int_0^t
e^{\lambda(t-s)}\,\CCO{R(\lambda,A)-R(\lambda,P_nA)}\,
P_ng(s)\,d\wien(s) \,d\lambda\\
&=\frac{\rho}{2\pi}\int_0^{2\pi}\int_0^t
e^{\rho\, e^{i\theta}\,(t-s)}\,\CCO{R(\rho\, e^{i\theta},A)-R(\rho\, e^{i\theta},P_nA)}\,
P_ng(s)\,d\wien(s) \,d\theta
\end{aligned}
\end{multline*}
by the stochastic Fubini theorem. 
Denoting again by $C$ a constant which may change from line to line,
we thus have 
\begin{multline*}
\expect\sup_{0\leq t\leq \Time}
\norm{
\int_0^t \CCO{e^{(t-s)A}-e^{(t-s)P_nA}}P_ng(s)\,d\wien(s)  
}^p\\
\begin{aligned}
&\leq C \expect\int_0^{2\pi}\sup_{0\leq t\leq \Time}
 \int_0^t\norm{
 e^{\rho\, e^{i\theta}\,(t-s)}\,\CCO{R(\rho\, e^{i\theta},A)-R(\rho\, e^{i\theta},P_nA)}\,
 P_ng(s)\,d\wien(s)
 }^p d\theta\\
&\leq C \expect\int_0^{2\pi}
 \sup_{0\leq t\leq \Time}\int_0^t\norm{
 e^{\rho\, e^{i\theta}\,(t-s)}\,\CCO{R(\rho\, e^{i\theta},A)-R(\rho\, e^{i\theta},P_nA)}\,
 P_ng(s)\,d\wien(s)
 }^p d\theta\\
&\leq C \int_0^{2\pi}\int_0^\Time
 \expect\norm{
 \CCO{R(\rho\, e^{i\theta},A)-R(\rho\, e^{i\theta},P_nA)}\,
 P_ng(s)\,d\wien(s)
 }^p d\theta\\
\end{aligned}
\end{multline*}
using the convolution inequality for the semigroup 
$t\mt e^{\rho\, e^{i\theta}\,t}$. 
From the compactness of $\comp$ and 
Lebesgue's
dominated convergence theorem, we get 
$$
\sup_{g\in\Gamma}\expect\sup_{0\leq t\leq \Time}
\norm{
\int_0^t \CCO{e^{(t-s)A}-e^{(t-s)P_nA}}P_ng(s)\,d\wien(s)  
}^p
$$
and we conclude as in Step 3.

\medskip
{\em Fifth step: Assuming all
preceding reductions, we now prove Lemma \ref{lem:Zero}.} 

Recall that $R_1=\sup_{x\in\comp}\norm{x}_\HSz$. 
We have, for any $\epsilon>0$ and $R>0$, 
\begin{align*}
\pr\biggl\{\sup_{0\leq t\leq\Time}\norm{\soll(t)}\leq R\biggr\}
&\leq \dfrac{4}{R^2}\expect\norm{\soll(t)}^2  \\
&=    \dfrac{4}{R^2}\expect\int_0^\Time \norm{e^{(\Time-s)A}g(s)}^2\,ds \\
&\leq \dfrac{4}{R^2}\Time e^{\Time\norm{A}} R_1^2  
\end{align*} 
Taking $R$ large enough, we get
\begin{equation}
\foreach{\epsilon>0}
\thereis{R>0}
\pr\biggl\{\sup_{0\leq t\leq\Time}\norm{\soll(t)}\leq R\biggr\} 
\leq \epsilon.\label{eq:aldousI} 
\end{equation}

Now, let $\tas$ be the set of stopping times $\tau$ such that
$0\leq\tau\leq\Time$. If $\sigma,\tau\in\tas$ with
$0<\tau-\sigma\leq\delta$ for some $\delta>0$, 
we have, for any $\soll\in\Soll$ of the form $\soll=\soll^{(g)}$,
with $g\in G$, and for any $\eta>0$, 
\begin{align*}
\pr\accol{\norm{\soll(\tau)-\soll(\sigma)}>\eta}
&\leq \dfrac{1}{\eta^2}\expect\norm{\soll(\tau)-\soll(\sigma)}^2\\
&\leq \dfrac{2}{\eta^2}
\begin{aligned}[t]
\biggl(
\expect\int_0^\sigma
\norm{\bigl(e^{(\tau-s)A}-e^{(\sigma-s)A}\bigr)g(s)}^2\,ds \phantom{\biggr)}&\\
\phantom{\leq \dfrac{2}{\eta^2}\biggl(}
+\expect\int_\sigma^\tau \norm{e^{(\tau-s)A}g(s)}^2\,ds
\biggr)&
\end{aligned}\\
&\leq \dfrac{2}{\eta^2}\CCO{\Time\delta e^{2\Time\norm{A}}\norm{A}^2
  R_1^2 
+\delta e^{\delta{\norm{A}}R_1^2}}\\
&\leq\delta\, C(\eta)
\end{align*}
{with} 
$C(\eta)=\dfrac{2e^{2\Time\norm{A}}\,R_1^2\,(1+\Time)}{\eta^2}$. 
Taking $\delta$ small enough, we get
\begin{equation}\label{eq:aldousII}
\foreach{\epsilon>0}
\foreach{\eta>0}
\thereis{\delta>0}
\foreach{\soll\in\Soll}
\sup_{\substack{\sigma,\tau\in\tas\\0< \tau-\sigma\leq\delta}}
\pr\accol{\norm{\soll(\tau)-\soll(\sigma)}>\eta}
\leq\epsilon. 
\end{equation}
From \eqref{eq:aldousI} and
\eqref{eq:aldousII}, we conclude, by a criterion of
Aldous \cite{aldous78stopping,jacod85limite}, 
that $\Soll$ is tight. 
\fin

\paragraph{The multivalued operator $\Phi$}
Let us denote by 
$\Phi$ the mapping which, with every continuous adapted $\esp$--valued process
$\sol$ such that $\expect{\int_0^\Time\norm{\sol(s)}^p\,ds}<+\infty$,  
associates the set of all processes of the form 
$$\semigr{t}\initial+\int_0^t
\semigr{t-s}f(s)\,ds+\int_0^t 
\semigr{t-s} g(s)\,d\wien(s),$$
where $f$ and $g$ are predictable selections of 
$(\omega,t)\mt \allowbreak F(t,\allowbreak \sol(\omega,t))$ and 
$(\omega,t)\mt \allowbreak G(t,\allowbreak \sol(\omega,t))$ respectively. 
Lemma \ref{lem:Zero} will be used through the following corollary: 
\begin{cor}\label{cor:Zero}
\titre{The operator $\Phi$ maps tight sets into tight sets}
Let $\EProc$ be a set of continuous adapted processes on
 $\esp$. Assume that each element of $\EProc$ is in
 $\Ellp{p}(\esprob\times\Times;\esp)$ and that $\EProc$, considered
 as a set of $\skoroC(\Times;\esp)$--valued random variables, is tight.  
 Assume furthermore that Hypothesis \HySg\ and \HyF\ are
 satisfied. 
 Then 
 $\Phi\circ\EProc:=\cup_{\proc\in\EProc}\Phi(\proc)$ 
 is a tight set of $\skoroC(\Times;\esp)$--valued random
 variables. 
\end{cor}
\preuve
Let us denote by $\Phi_F$ the mapping which, 
with every continuous adapted $\esp$--valued process
$\sol$ such that $\expect{\int_0^\Time\norm{\sol(s)}^p\,ds}<+\infty$,  
associates the set of all processes of the form 
$$
\int_0^t \semigr{t-s}f(s)\,ds,
$$
where $f$ is a predictable selection of 
$(\omega,t)\mt \allowbreak F(t,\allowbreak \sol(t))$. 

By Ascoli's theorem (see the details in the proof of \cite[Lemma 4.2.1]{koz01book}), 
the mapping 
$$ f\mt \int_0^{\displaystyle .}
\semigr{.-s}f(s)\,ds $$
maps all measurable functions 
$f :\,\Times\ra \esp$ with values in a given compact subset of
$\esp$ into a compact subset of $\skoroC(\Times;\esp)$.  
Thus, following the same lines as in Steps 1 and 2 of the proof of
Lemma \ref{lem:Zero}, $\Phi_F$
maps tight bounded subsets of
$\Ellp{p}(\esprob\times\Times;\esp)$ into tight sets of
$\skoroC(\Times;\esp)$--valued random  variables. 
The conclusion immediately follows from Lemma
 \ref{lem:Zero}. 
\fin

\begin{lem}\label{lem:Aprime}
\titre{The operator $\Phi$ and the measure of noncompactness $\mncp$}
Let $\EProc$ be a bounded subset of $\BTC{\Time}{\esp}$. Assume 
Hypothesis 
\HyFG. 
We then have 
\begin{equation*}
{\mncp}^p(\Phi\circ\EProc)(t)
\leq \Lcte  
\int_0^t \LL\CCO{s,{\mncp}^p(\EProc)}(s)\,ds
\end{equation*}
for some constant $\Lcte$ which depends only on 
$\Time$, $p$, $\CteCauchy{\Time}$,and $\CteConv$. 
\end{lem}
\preuve
The main arguments are inspired from \cite[Lemma 4.2.6]{akprs92book}. 
For simplicity, the space $\BTC{t}{\esp}$ will be denoted by
$\BTCsimple{t}$. 

Let $\epsilon>0$. 
The function $t\mt\mncp(\EProc)(t)$ is increasing, thus there exist at
most a finite number of points $0\leq t_1\leq\dots\leq t_n\leq\Time$
for which $\mncp(\EProc)$ makes a jump greater than $\epsilon$. 
Let us choose $\delta_1>0$ such that 
$i\not=j\Rightarrow ]t_i-\delta_1,t_i+\delta_1[\cap
]t_j-\delta_1,t_j+\delta_1[=\emptyset$. 
Using points $\beta_j$, $j=1,\dots,m$, 
we divide the remaining part 
 $[0,\Time]\setminus \CCO{
]t_1-\delta_1,t_1+\delta_1[\cup\dots\cup]t_n-\delta_1,t_n+\delta_1[}$  
into disjoint intervals in such
a way that, for each $j$, 
\begin{equation}
  \label{eq:diameter-range}
  \sup_{s,t\in[\beta_{j-1},\beta_j]}\abs{\mncp(\EProc)(s)-\mncp(\EProc)(t)}<\epsilon. 
\end{equation}
Let us then choose $\delta_2>0$ such that 
$j\not=k\Rightarrow ]\beta_j-\delta_2,\beta_j+\delta_2[\cap
]\beta_k-\delta_2,\beta_k+\delta_2[=\emptyset$.

Now we start the construction of a tight net. 

For each $j=1,\dots,m$, we choose a tight 
$\CCO{{\mncp}(\EProc)(\beta_j)+\epsilon}$--net $N_j$ of $\EProc$
in $\BTCsimple{\beta_j}$.  
(As $\BTCsimple{\beta_j}$ is separable, we can take a countable net
$N_j$, but this fact will not be used here.)  
We obtain a (countable) family $\ZZ$ in $\BTCsimple{\Time}$ by
taking all continuous processes which coincide on each
$]\beta_{j-1}+\delta_2,\beta_j-\delta_2[$ ($1\leq j\leq m$) with some
element of $N_j$ and which have affine trajectories on the
complementary segments. 
The set $\ZZ$ is tight and bounded in
$\Ellp{p}(\esprob\times\Times;\esp)$, 
thus, by Corollary \ref{cor:Zero},  
$\Phi(\ZZ)$ is tight. 

Consider a fixed $\sol\in\EProc$. 

We can find an element $Z$ of $\ZZ$ such that, for each $j=1,\dots,m$, 
\begin{equation*} 
\norm{\sol-Z}_{\BTCsimple{\beta_j}}
\leq
{{\mncp}(\EProc)(\beta_j)+\epsilon}. 
\end{equation*}
For $t\in]\beta_{j-1}+\delta_2,\beta_j-\delta_2[$, we have, using
\eqref{eq:diameter-range},  
\begin{align}
\expect\norm{\sol(t)-Z(t)}_\esp^p
&\leq \expect \sup_{\beta_{j-1}+\delta_2\leq s\leq \beta_j+\delta_2}
\norm{\sol(s)-Z(s)}_\esp^p \notag \\
&\leq \norm{\sol-Z}_{\BTCsimple{\beta_j}}^p\notag \\
&\leq \CCO{{\mncp}(\EProc)(\beta_j)+\epsilon}^p\notag \\
&\leq \CCO{{\mncp}(\EProc)(t)+2\epsilon}^p \label{eq:ptitsurI}
\end{align}

Let $\soll\in\Phi(\sol)$, say 
$$\soll(t)=\semigr{t}\sol(0)+\int_0^t
\semigr{t-s}f(s)\,ds+\int_0^t 
\semigr{t-s} g(s)\,d\wien(s),$$
where $f$ and $g$ are predictable selections of 
$(\omega,t)\mt \allowbreak F(t,\allowbreak \sol(\omega,t))$ and 
$(\omega,t)\mt \allowbreak G(t,\allowbreak \sol(\omega,t))$
respectively. 
We can find predictable selections 
$\tilde{f}$ and $\tilde{g}$ of  $(\omega,t)\mt \allowbreak F(t,\allowbreak Z(\omega,t))$ and 
$(\omega,t)\mt \allowbreak G(t,\allowbreak Z(\omega,t))$
respectively, such that, for every $t\in\Times$, 
\begin{align}
\norm{\tilde{f}(t)-f(t)}^p
&\leq 2\Hausd^p\bigl(F(t,\sol(t)),F(t,Z(t))\bigr)\notag\\
&\leq 2\LL\bigl(t,\norm{\sol(t)-Z(t)}^p\bigr)
\label{eq:fftildeL}\\
\intertext{and}
\norm{\tilde{g}(t)-g(t)}^p
&\leq 2\Hausd^p\bigl(G(t,\sol(t)),G(t,Z(t))\bigr)\notag\\
&\leq 2\LL\bigl(t,\norm{\sol(t)-Z(t)}^p\bigr). 
\label{eq:ggtildeL}
\end{align} 
Let $d_{\BTCsimple{t}}$ be the distance in $\BTCsimple{t}$ associated
with $\norm{.}_{\BTCsimple{t}}$. 
We have, using \HyFG\ and the convexity of 
$x\mt x^p$, 
\begin{align*}
d^p_{\BTCsimple{t}}\CCO{\soll,\Phi(Z)}
&\leq \expect \sup_{0\leq \tau\leq t}\,
\bigl\Vert
\int_0^\tau \semigr{t-s}(f(s)-\tilde{f}(s))\,ds\\
&\phantom{\leq \expect \sup_{0\leq \tau\leq t}\,\bigl\Vert}
+\int_0^\tau \semigr{t-s}(g(s)-\tilde{g}(s))\,d\wien(s)
\bigr\Vert^p\\
&\leq 2^{p-1}\expect\,\sup_{0\leq \tau\leq t}\,
\bigl\Vert \int_0^\tau \semigr{t-s}(f(s)-\tilde{f}(s))\,ds
\bigr\Vert^p\\
&\phantom{\leq \expect \sup_{0\leq \tau\leq t}\,\bigl\Vert}
+\expect\,\sup_{0\leq \tau\leq t}\,
\Vert \int_0^\tau \semigr{t-s}(g(s)-\tilde{g}(s))\,d\wien(s)
\bigr\Vert^p\\
&\leq 2^{p-1}\bigl(
\CteCauchy{t}\int_0^t\expect \norm{f(s)-\tilde{f}(s)}^p\,ds\\
&\phantom{\leq \expect \sup_{0\leq \tau\leq t}\,\bigl\Vert}
+\CteConv\,\Time^{p/2-1}\int_0^t\expect \norm{g(s)-\tilde{g}(s)}^p\,ds
\bigr)
\end{align*}
(recall that $\CteCauchy{t}$ is defined in \HySg). 
{Let us denote}
\begin{gather}
\Lcte=2^{p}\max\CCO{\CteCauchy{\Time},\CteConv\,\Time^{p/2-1}},\label{eq:def-de-k}\\
J(t)=[0,t]\cap
 \biggl((\cup_{1\leq i\leq n}]t_i-\delta_1,t_i+\delta_1[)
 \cup(\cup_{1\leq j\leq
   m}]\beta_j-\delta_2,\beta_j+\delta_2[)\biggr),\notag \\
I(t)=[0,t]\setminus J(t).\notag
\end{gather}
As $\EProc$ and $\ZZ$ are bounded in
$\Ellp{p}(\esprob\times\Times;\esp)$, using \HyFG\ 
\eqref{eq:fftildeL} and \eqref{eq:ggtildeL}, 
and taking
$\delta_1$ and $\delta_2$ sufficiently small, we get 
\begin{align*}
d^p_{\BTCsimple{t}}\CCO{\soll,\Phi(Z)}
&\leq
\Lcte
\bigl(\int_{I(t)} \expect \LL(s,\norm{\sol(s)-Z(s)}^p)\,ds 
+\epsilon \bigr).
\end{align*}
Then, using the convexity of 
$\LL(t,.)$, and \eqref{eq:ptitsurI},
\begin{align*}
d^p_{\BTCsimple{t}}\CCO{\soll,\Phi(Z)}
 &\leq \Lcte\bigl( \int_0^t\LL\CCO{s,\CCO{{\mncp}(\EProc)(s)+2\epsilon}^p}\,ds
 +\epsilon\bigr).
\end{align*}
As $\epsilon$, $\sol$ and $\soll$ are arbitrary, the result follows. 
\fin
Here is a easy variant for the case when $F$ and $G$ are single-valued:
\begin{lem}\label{lem:compact-single}
Assume that $F$ and $G$ are single-valued. Assume furthermore Hypothesis 
\HyFG. 
Let $\mncpp$ be the measure of noncompactness on $\BTC{\Time}{\esp}$ 
defined by 
$$   \mncpp(\EProc)(s):=\inf\accol{\epsilon>0\tq
   \text{$\EProc\restr{s}$ has a finite $\epsilon$--net in
   $\BTC{s}{\esp}$}}
   \quad (0\leq s\leq\Time).$$
Let $\EProc$ be a bounded subset of 
$\BTC{\Time}{\esp}$. 
We have 
\begin{equation*}
{\mncpp}^p(\Phi\circ\EProc)(t)
\leq \Lcte'  
\int_0^t \LL\CCO{s,{\mncpp}^p(\EProc)}(s)\,ds
\end{equation*}
for some constant $\Lcte'$ which depends only on 
$\Time$, $p$, $\CteCauchy{\Time}$,and $\CteConv$. 
\end{lem}
\preuve 
We only need to repeat the proof of Lemma \ref{lem:Aprime}, but we take for
each $j$ a finite $\CCO{{\mncpp}(\EProc)(\beta_j)+\epsilon}$--net
$N_j$. Then $\ZZ$ is finite, thus $\Phi(\ZZ)$ is compact in
$\BTC{\Time}{\esp}$. The rest of the proof goes as in the proof of 
Lemma \ref{lem:Aprime}. 
\fin

In the proof of Theorem \ref{theo:main}, we shall consider a variant
of the operator $\Phi$. 
For each $n\geq 1$, let $\Phi_n$ 
be the mapping which, with every continuous adapted $\esp$--valued process
$\sol$ such that $\expect{\int_0^\Time\norm{\sol(s)}^p\,ds}<+\infty$, 
associates the set of all processes of the form 
$$\semigr{t-1/n}\initial+\int_0^{t-1/n}
\semigr{t-s}f(s)\,ds+\int_0^{t-1/n} 
\semigr{t-s} g(s)\,d\wien(s),$$
where $f$ and $g$ are predictable selections of 
$(\omega,t)\mt F(t,\sol(\omega,t))$ and 
$(\omega,t)\mt G(t,\sol(\omega,t))$. 
The following lemma links the tightness properties of $\Phi$ and
$\Phi_n$. 
\begin{lem}\label{lem:C}
Let $(\sol_n)$ be a sequence of continuous adapted $\esp$--valued
processes, which is bounded in $\Ellp{p}(\esprob\times\Times;\esp)$. 
We then have 
\begin{equation*}
\mncp\CCO{\cup_{n}\Phi_n(\sol_n)}
\leq \CteCauchy{\Time}\mncp\CCO{\cup_n\Phi(\sol_n)},
\end{equation*}
where $\CteCauchy{\Time}$ has been defined in \HySg. 
\end{lem}
\preuve
First, in the definition of $\Phi_n$ and in that of $\Phi$, we can
assume that $\initial=0$, because $\initial$ does not change the
values of $\mncp\CCO{\cup_{n}\Phi_n(\sol_n)}$ and
$\mncp\CCO{\cup_n\Phi(\sol_n)}$. 
We thus have 
\begin{align*}
\Phi_n(\sol_n)(t)
&=\biggr\{\begin{aligned}[t]
\semigr{{1}/{n}}\int_0^{t-1/n}&\semigr{(t-1/n)-s}f(s)\,ds\\
&+\semigr{{1}/{n}}\int_0^{t-1/n}\semigr{(t-1/n)-s}g(s)\,d\wien(s)
\tq \\
&\text{ $f$ and $g$ are predictable selections of $F\circ\sol$ and $G\circ\sol$}
\biggl\}\end{aligned}\\
&=\semigr{{1}/{n}}\Phi(\sol_n)(t-1/n). 
\end{align*}
Let $\tau_n :,\skoroC(\Times;\esp)\ra\skoroC(\Times;\esp)$ be defined by 
$$\tau_n(u)(t)=
\begin{cases}
u(t-1/n)& \text{ if }t\geq 1/n\\
u(0)    & \text{ if }t\leq 1/n. 
\end{cases}
$$
For any set $\EProc$ of continuous processes, we have 
$$\mncp(\tau_n(\EProc))\leq \mncp(\EProc)$$ 
because, if $\Xi$ is a tight $\epsilon$--net of $\EProc$, 
then $\tau_n(\Xi)$ is a tight $\epsilon$--net of $\tau_n(\EProc)$. 
We thus have 
\begin{multline*}
\rule{1cm}{0em}
\mncp\CCO{\cup_{n}\Phi_n(\sol_n)}
=\mncp\CCO{\semigr{{1}/{n}}\cup_n\tau_n(\Phi(\sol_n))}\\
\leq \CteCauchy{\Time} \mncp\CCO{\cup_n\tau_n(\Phi(\sol_n))}
\leq \CteCauchy{\Time} \mncp\CCO{\cup_n\Phi(\sol_n)}. 
\rule{1cm}{0em}
\end{multline*}
\fin

\section{Proof of the main results} 
\label{sect:proofs}

\preuvof{Theorem \ref{theo:main}}

{\em First step: construction of a tight sequence of
approximating solutions through Tonelli's scheme}. 


For each integer $n\geq 1$, we can easily define a
process $\ssol_n$ on [-1,\Time] by 
$\ssol_n(t)=0$ if $t\leq 0$
and, for $t\geq 0$,  
\begin{align*}
   \ssol_n(t)
  &=  \semigr{t}\initial+\int_0^t\semigr{t-(s-\frac{1}{n})}f_n(s)\,ds 
     +\int_0^t\semigr{t-(s-\frac{1}{n})}g_n(s)\,d\wien(s), 
\end{align*}
where $f_n :\,\esprob\times\Times\times\ra\esp$ and 
      $g_n :\,\esprob\times\Times\times\ra\HSz$ are predictable and 
\begin{align*}
  f_n(s)&\in F(s,\ssol_n(s-1/n)) \text{ a.e.}\quad\text{ and }\quad
  g_n(s)\in G(s,\ssol_n(s-1/n)) \text{ a.e.}\
\end{align*}
We then set, for $t\in[1/n,\Time]$,  
\begin{align*}
  \sol_n(t)
  &= \ssol_n(t-1/n)\\
  &= \semigr{t-1/n}\initial+\int_0^{t-1/n}\semigr{t-s}f_n(s)\,ds 
     +\int_0^{t-1/n}\semigr{t-s}g_n(s)\,d\wien(s). 
\end{align*}
For $t\leq 1/n$, we set $\sol_n(t)=\initial$. 
Let us show that $(\sol_n)$ is bounded in
$\Ellp{p}(\esprob\times\Times;\esp)$. 
By convexity of $t\mt\abs{t}^p$, we have 
\begin{equation}
\norm{a_1+\dots+a_m}^p\leq m^{p-1}\CCO{\norm{a_1}^p+\dots+\norm{a_m}^p}
\end{equation}
for any finite sequence $a_1,\dots,a_m$ in any normed space. 
Recall also that $\CteCauchy{\Time}$ and $\growthconstant$ have been
  defined in Hypothesis \HySg\ and \HyF-\qref{hyp:growth} and
  $\CteConv$ is the constant of stochastic convolution defined in
  \eqref{stochconvcondense}.  
For every $n\geq 1$, we have the following chain of inequalities,
where the supremum is taken over  all predictable
  selections $(f,g)$ of $(\omega,t)\mt F(t,\sol_n(t))\times
  G(t,\sol_n(t))$:
\begin{align*}
\expect{\norm{\sol_n(t)}^p}
&\leq \sup_{f,g}
3^{p-1}\begin{aligned}[t]\biggl(
      \expect\norm{\semigr{t-1/n}\initial}
      &+\expect\int_0^{t-1/n}\norm{\semigr{t-s}f(s)}^p\,ds\\
      &+\expect\int_0^{t-1/n}\norm{\semigr{t-s}g(s)}^p\,d\wien(s)
           \biggr)\end{aligned}\\
&\leq \sup_g 3^{p-1}\begin{aligned}[t]\biggl(
      \CteCauchy{\Time}^p\expect\norm{\initial}^p
      &+\CteCauchy{\Time}^p\growthconstant^p\expect\int_0^{t}\CCO{1+\norm{\sol_n}(s)}^p\,ds\\
      &+\CteCauchy{\Time}^p\CteConv\expect\int_0^{t-1/n}\norm{g(s)}^p\,ds
           \biggr)\end{aligned}\\
&\leq 3^{p-1}\CteCauchy{\Time}^p\begin{aligned}[t]\biggl(
       \expect\norm{\initial}^p
       &+\growthconstant^p(1+\CteConv)\expect\int_0^t\CCO{1+\norm{\sol_n(s)}}^p\,ds
          \biggr).\end{aligned}\\
&\leq 3^{p-1}\CteCauchy{\Time}^p\begin{aligned}[t]\biggl(
       \expect\norm{\initial}^p
       &+\growthconstant^p(1+\CteConv)2^{p-1}\expect\int_0^t 1+\norm{\sol_n(s)}^p\,ds
          \biggr).\end{aligned}
\end{align*}
Let
$\grosscte
=3^{p-1}\CteCauchy{\Time}^p
\bigl(\expect\norm{\initial}^p+2^{p-1}\growthconstant^p(1+\CteConv)\Time\bigr)$.  
We have 
$$
\expect{\norm{\sol_n(t)}^p}
\leq\grosscte +\grosscte\int_0^t\expect\norm{\sol_n(s)}^p\,ds,
$$
thus, by Gronwall Lemma,
$$
\expect{\norm{\sol_n(t)}^p}\leq \grosscte e^{\grosscte t}$$
which provides the boundedness condition 
\begin{equation}\label{eq:xn-bounded}
\sup_n\expect\int_0^\Time \norm{\sol_n(t)}^p\,dt
\leq \Time\grosscte e^{\grosscte \Time}. 
\end{equation}

Let us now show that $(\sol_n)$ is a tight sequence of
$\skoroC(\Times;\esp)$--valued random variables. 
As $(\sol_n)$ is bounded in
$\Ellp{p}(\esprob\times\Times;\esp)$, we can apply Lemma \ref{lem:C}
and Lemma \ref{lem:Aprime}. 
Using  
the fact that $\sol_n\in\Phi_n(\sol_n)$ for each $n$ , we get, for every
$t\in\Times$,  
\begin{align*}
\CCO{\mncp(\cup_n\{\sol_n\})(t)}^p
&\leq \CCO{\mncp(\cup_n\Phi_n(\sol_n)(t)}^p
\leq \CteCauchy{\Time}^p\mncp\CCO{\cup_n\Phi(\sol_n)(t)}^p\\
&\leq \CteCauchy{\Time}^p\,\Lcte  
\int_0^t \LL\CCO{s,{\mncp}^p(\cup_n\{\sol_n\})(s)}\,ds
\end{align*}
(where $\Lcte$ is the constant we obtained in Lemma \ref{lem:Aprime},
see \eqref{eq:def-de-k}). 
Thus, by \eqref{eq:hyp:L} in Hypothesis \HyF-\qref{hyp:L}, 
we have 
${\mncp(\cup_n\{\sol_n\})(t)}=0$ for each $t$, which, by Lemma
\ref{lem:G}, implies that $(\sol_n)$ is tight.

\medskip
{\em Second step: construction of a weak solution}. 

By Prohorov's compactness criterion for Young measures, 
we can extract a subsequence of $(X_n)$ which converges 
{stably} 
to a Young measure $\mu\in\youngs(\esprob,\tribu,\pr;\skoroC(\Times;\esp))$. 
For simplicity, we denote this extracted sequence by $(X_n)$.

It will be convenient to represent the limiting Young measure $\mu$ as a
random variable defined on an extended probability space. 
Let $\Ctribu$ be the Borel $\sigma$-algebra of $\skoroC(\Times;\esp)$
and, for each $t\in\Times$, let 
$\Ctribu_t$ be the sub-$\sigma$--algebra of $\Ctribu$ generated by 
$\skoroC([0,t];\esp)$. 

We define a stochastic basis $(\esprobb,\tribuu,(\tribuu_t)_t,\mu)$ by 
\begin{equation*}
\esprobb=\esprob\times\skoroC(\Times;\esp),
\quad 
\tribuu=\tribu\otimes\Ctribu
\quad 
\tribuu_t=\tribu_t\otimes\Ctribu_t
\end{equation*}
and we define 
$X_\infty$ on $\esprobb$ by 
$$X_\infty(\omega,u)=u.$$
Clearly, $\law{X_\infty}=\mu$ and $X_\infty$ is 
$(\tribuu_t)$--adapted. 
Now, the random variables $X_n$ can be seen as random elements defined
on $\esprobb$, using the notation
$$X_n(\omega,u):=X_n(\omega)\quad(n\in\N).$$
Furthermore, $X_n$ is $(\tribuu_t)$--adapted for each $n$. 
The $\sigma$--algebra $\tribu$ can be identified with the
sub--$\sigma$--algebra 
$\{A\times\skoroC(\Times;\esp)\tq A\in\tribu\}$ 
of $\tribuu$. 
We thus have:
\begin{equation}\label{eq:F-stable}
\foreach{A\in\tribu}
\foreach{f\in\Cb{\esp}}
\lim_n\int_A f(X_n)\,d\mu
=\int_A f(X_\infty)\,d\mu.
\end{equation}
To express \eqref{eq:F-stable}, we say that 
$(X_n)$ converges to $X_\infty$ {\em $\tribu$--stably}. 

Let us show that $\wien$ is an $(\tribuu_t)$--Wiener process under the
probability $\mu$. 
Clearly, $\wien$ is $(\tribuu_t)$--adapted,  
so we only need to prove
that $\wien$ has independent increments. 
By a result of Balder \cite{balder89Prohorov,balder90new}, 
each subsequence of $(\ydirac{X_n})$ contains a further subsequence 
$(\ydirac{{X}'_n})$
which K--converges to $\mu$, that is, for each subsequence 
$(\ydirac{{X}''_n})$ of $(\ydirac{{X}'_n})$, we have 
$$\lim_n\frac{1}{n}\sum_{i=1}^n\dirac{{X}''_n(\omega)}=\mu_\omega\ \text{a.e.}$$
This entails that, for every $A\in\Ctribu_t$, the mapping $\omega\mt\mu_\omega(A)$
is $\tribu_t$--measurable\footnotemark. 
\footnotetext{From \cite[Lemma 2.17]{jacod-memin81weakstrong}, 
this means that $(\esprobb,\tribuu,(\tribuu_t)_t,\mu)$ 
is a {\em very good extension} 
of $(\esprob,\tribu,(\tribu_t)_t,\pr)$
in the sense of
\cite{jacod-memin81weakstrong}, that is, 
every martingale on $(\esprob,\tribu,(\tribu_t)_t,\pr)$ remains a martingale 
on $(\esprobb,\tribuu,(\tribuu_t)_t,\mu)$.} 
Let $t\in\Times$ and let $s>0$ such that $t+s\in\Times$. 
Let us prove that, for 
any $A\in\tribuu_t$ and any  Borel subset $C$ of $\hilbw$, we have 
\begin{equation}\label{eq:independence}
\mu\CCO{A\cap \{\omega\in\esprob\tq \wien(t+s)-\wien(t)\in C\}}
=\mu(A)\, \mu\{\omega\in\esprob\tq \wien(t+s)-\wien(t)\in C\}.
\end{equation} 
Let $B=\{\omega\in\esprob\tq \wien(t+s)-\wien(t)\in C\}$. 
We have
\begin{align*}
\mu\CCO{A\cap (B\times\skoroC(\Times;\esp))}
&=\int_{\Omega\times\skoroC(\Times;\esp)} 
\un{A}(\omega,u))\un{B}(\omega)\,d\mu(\omega,u) \\
&=\int_\Omega\mu_\omega(\un{A}(\omega,.))\un{B}(\omega)\,d\pr(\omega) \\
&=\int_\Omega\mu_\omega(\un{A}(\omega,.))\,d\pr(\omega)\,\pr(B) \\
&=\mu(A)\, \mu(B\times\skoroC(\Times;\esp)),
\end{align*} 
which proves \eqref{eq:independence}. Thus 
$\wien(t+s)-\wien(t)$ is independent of $\tribuu_t$. 

Now, there remains to prove that $X=X_\infty$ satisfies \eqref{eq:generale}. 
Note that the first step of the proof of Theorem \ref{theo:main} 
is valid for any choice of the selections $f_n$ and
$g_n$. In this second part, we will use a particular choice. 
As $F$ and $G$ are globally measurable, they have measurable graphs 
(see \cite[Proposition III.13]{castaingvaladier77book}). As
furthermore they are    
continuous in the second 
variable, they admit Carath\'eodory selections 
\cite[Corollary 1 of the Main Lemma]{kucia98carath},  
that is, there exist globally
measurable mappings 
$\carathf:\,\Times\times\esp\rightarrow \esp $ and 
$\carathg :\,\Times\times\esp\rightarrow \HSz $ 
such that $\carathf(t,.)$ and $\carathg(t,.)$ are continuous for every
$t\in\Times$ and $\carathf(t,x)\in F(t,x)$ and $\carathg(t,x)\in
G(t,x)$ for all 
$(t,x)\in\Times\times\esp$. 
We denote by $\NN$ the one-point compactification of $\N$ and we set,
for every $n\in\NN$ and every 
$s\in\Times$, 
\begin{alignat*}{2}
f_n(s)&=\carathf(s,X_n(s))\\
g_n(s)&=\carathg(s,X_n(s)). 
\end{alignat*}
The sequence $(f_n)$ converges in law to $f_\infty$ in the space 
$\Ellp{p}\CCO{\Times;\esp}$. 
Indeed, 
by the growth condition \HyF-\qref{hyp:growth} and Lebesgue's
convergence theorem, 
the mapping 
$$\supercarath:\,\left\{\begin{array}{lcl}
\skoroC(\Times;\esp])&\rightarrow\Ellp{p}\CCO{\Times;\esp}\\
u&\mapsto \carathf(.,(u(.)))
\end{array}\right.$$
is well defined and continuous, and we have 
$f_n=\supercarath\circ X_n$. 
Similarly, $(g_n)$ converges in law to 
$g_\infty$ in $\Ellp{p}\CCO{\Times;\HSz}$. 
Actually, we have more: the sequence 
$(X_n,\wien,f_n,g_n)$ converges in law to 
$(X_\infty,\wien,f_\infty,g_\infty)$ in 
$\skoroC(\Times;\esp)$\allowbreak$\times$\allowbreak
$\skoroC(\Times;\hilbw)$\allowbreak$\times$\allowbreak
$\Ellp{p}\CCO{\Times;\esp}$\allowbreak$\times$\allowbreak$\Ellp{p}\CCO{\Times;\HSz}$. 
In particular, 
as $p>2$, the convergence in law of $(\wien,g_n)$ to $(\wien,g_\infty)$
implies, for every $t\in\Times$, the convergence in law of the stochastic
integrals 
$
\int_0^t\semigr{t-s}g_n(s)\,d\wien(s)
$
to 
$
\int_0^t\semigr{t-s}g_\infty(s)\,d\wien(s). 
$
To see this, one can e.g.~use a Skorokhod representation $(\wien'_n,g'_n)$ which
converges a.e.~to $(\wien',g_\infty')$, then, by Vitali's convergence
theorem,  $(\wien'n,g'_n)$ converges in
$\Ellp{2}\CCO{\esprob\times\Times;\hilbw\times\HSz}$ thus, for every
$t\in\Times$ the stochastic integrals
$
\int_0^.\semigr{t-s}g'_n(s)\,d\wien'(s)
$
converge in 
$\Ellp{2}\CCO{\esprob\times[0,t];\esp}$.  
Let $t\in\Times$ be fixed and 
let us set, for every $n\in\NN$,  
$$Z_n(t)=-X_n(t)+\semigr{t}\initial
+\int_0^{t-1/n}\semigr{t-s}f_n(s)\,ds
+\int_0^{t-1/n}\semigr{t-s}g_n(s)\,d\wien(s)$$
(with $1/\infty:=0$). 
For every $t\in\Times$, the sequence $(Z_n(t))$ converges in law
to $Z_\infty(t)$. 
Now, from the definition of $X_n$ ($n<+\infty$), we have 
\begin{multline*}
\foreach{n<+\infty}
Z_n(t)=\CCO{\semigr{t}-\semigr{t-1/n}}\initial\\
+\int_{t-1/n}^t\semigr{t-s}f_n(s)\,ds
+\int_{t-1/n}^t\semigr{t-s}g_n(s)\,d\wien(s).
\end{multline*}
But $\CCO{\semigr{t}-\semigr{t-1/n}}\initial$ converges a.e.~to 0 and
we have, using  the growth condition
  \HyF-\qref{hyp:growth} and the boundedness property \eqref{eq:xn-bounded},
\begin{align*}
\expect\norm{\int_{t-1/n}^t\semigr{t-s}f_n(s)\,ds}
&\leq \sgcontract e^{\sgcoef \Time}\int_{0}^t\norm{f_n(s)\un{[t-1/n,t]}(s)}\,ds\\
&\leq \sgcontract e^{\sgcoef \Time}
\expect\CCO{\int_{0}^t\norm{f_n(s)}^p\,ds}^{1/p}\CCO{\frac{1}{n}}^{p/(p-1)}\\
&\leq \sgcontract e^{\sgcoef
  \Time}\CCO{2^{p-1}\expect\int_0^t(1+\norm{X_n(s)})^p\,ds}^{1/p}\CCO{\frac{1}{n}}^{p/(p-1)}\\
&\rightarrow 0 \text{ as }n\rightarrow \infty
\end{align*}
and, for any $q$ such that $2<q<p$,
\begin{multline*}
\expect\norm{\int_{t-1/n}^t\semigr{t-s}g_n(s)\,d\wien(s)}^q\\
\begin{aligned}
&=\expect\norm{\int_{0}^t\semigr{t-s}g_n(s)\un{[t-1/n,t]}(s)\,d\wien(s)}^q\\
&\leq \CteConv'\Time^{q/2}\expect\int_0^t
\norm{g_n(s)\un{[t-1/n,t]}(s)}q\,ds\\
&\leq \CteConv'\Time^{q/2}\CCO{\expect\int_0^t
  \norm{g_n(s)}^p\,ds}^{q/p}\CCO{\frac{1}{n}}^{(p-q)/p}\\
&\leq \CteConv'\Time^{q/2}\CCO{2^{q-1}\expect\int_0^t
  (1+\norm{g_n(s)}^p)\,ds}^{q/p}\CCO{\frac{1}{n}}^{(p-q)/p} \\
&\rightarrow 0 \text{ as }n\rightarrow \infty,
\end{aligned}
\end{multline*}
where $\CteConv'$ is the constant of stochastic convolution associated
with $q$. 
Thus $Z_n(t)$ converges to 0 in probability. 
Thus, for every $t\in\Times$, we have $Z_\infty(t)=0$ a.s. 
As $Z_\infty$ is continuous, this means that $Z_\infty=0$ a.s.  
Thus $X_\infty$ is a weak mild solution to \eqref{eq:generale}. 
\fin

\preuvof{Proposition \ref{prop:single-valued}}
Taking into account Lemma \ref{lem:compact-single}, 
and following the reasoning of the first part of the proof of Theorem
\ref{theo:main}, we see that the sequence $(X_n)$ provided by the
Tonelli scheme is relatively compact in $\BTC{\Time}{\esp}$. 
Moreover, the limit of any convergent  
subsequence of $(X_n)$ is a strong mild solution to \eqref{eq:generale}. 
\fin

If $\espgene$ is a Banach space,  
let us call {\em selection of $\compacts(\espgene)$} every mapping 
$\steiner :\,\compacts(\espgene)\rightarrow\espgene$ such that $\steiner(K)\in
K$ for every $K\in\compacts(\espgene)$. 

\preuvof{Proposition \ref{cor:finitedim}}
If $\HSz$ is finite dimensional, it is well-known 
that there exists a selection $\steiner_\HSz$ of $\compacts(\HSz)$
which is Lipschitz with respect to the Hausdorff distance. 
One such mapping is the Steiner
point (see 
e.g.~\cite{schneider71steiner,saint-pierre85steiner,vitale85steiner,przeslawski96centres}).  
In this case, we can define the Carath\'eodory mapping $\carathg$ of
the proof of Theorem \ref{theo:main} by 
$\carathg(s,x)=\steiner_\HSz(G(s,x))$. 
Then the selection $\carathg$ satisfies all hypothesis satisfied by $G$,
in particular \HyF-\qref{hyp:FL}.
Similarly, if $\esp$ is finite dimensional, we can define  
$\carathf$ as the Steiner point of $F$. 
So, the proof of Corollary \ref{cor:finitedim} 
reduces to the case 
when $F$ and $G$ are single valued, and, in this case, from Proposition
\ref{prop:single-valued}, there exists a strong mild solution to
\eqref{eq:generale}. 
\fin

  \begin{rem}   
  It is well known that, if $\HSz$ is infinite dimensional,
  there exists no Lipschitz selection  of $\compacts(\HSz)$: see   
  \cite[Theorem 4]{posicelski71lipschitzian}, 
  where the reasoning given for the set of convex
  bounded sets also applies to the set of convex compact sets, 
  and see also \cite{przeslawski-yost89continuity}. 
  \end{rem}




\medskip
\noindent{\bf Aknowledgements }
We thank Professor Jan Seidler for several constructive remarks 
and for making us know the papers 
\cite{yamada81successive,manthey90convergence,taniguchi92successive,barbu98local-global,barbu-bocsan02approx}.  

\def\cprime{$'$} \def\cprime{$'$}
  \def\polhk#1{\setbox0=\hbox{#1}{\ooalign{\hidewidth
  \lower1.5ex\hbox{`}\hidewidth\crcr\unhbox0}}}
  \def\polhk#1{\setbox0=\hbox{#1}{\ooalign{\hidewidth
  \lower1.5ex\hbox{`}\hidewidth\crcr\unhbox0}}} \def\cprime{$'$}
  \def\cprime{$'$} \def\cprime{$'$} \def\cprime{$'$} \def\cprime{$'$}
  \def\cprime{$'$}

\end{document}